\documentclass[11pt]{article}
\usepackage{mathtools}

\usepackage{amsmath,amssymb,amstext,amsfonts}
\usepackage{amsthm}
\usepackage{enumerate}
\usepackage{tcolorbox}
\usepackage{graphicx}
\usepackage{subcaption}
\usepackage{cite}
\usepackage{float}
\usepackage[margin=1 in]{geometry}
\usepackage{listings}
\usepackage{xcolor}
\usepackage{subcaption}

\usepackage{graphicx}
\usepackage[utf8]{inputenc}

\usepackage{systeme}
\usepackage{comment}
\usepackage{dsfont}
\newtheorem*{remark}{Remark}

\tcbuselibrary{theorems}
\theoremstyle{definition}

\newcounter{texercise}
\newwrite\solout
\def\openoutsol{\immediate\openout\solout\jobname.sol}  \def\writesol#1{\immediate\write\solout{\noexpand\processsol{\thetexercise}{#1}}}% \def\closeoutsol{\immediate\closeout\solout} \def\inputsol{\IfFileExists{\jobname.sol}{\input{\jobname.sol}}{}}
\tcbmaketheorem{theo}{Theorem}{colback=green!5,colframe=green!35!black, fonttitle=\bfseries}{texercise}{th}

\newcounter{mytheorem}[section] 
\tcbset{
theostyle/.style={fonttitle=\bfseries\upshape, fontupper=\slshape,
arc=0mm, colback=blue!5,colframe=blue!75!black}, defstyle/.style={fonttitle=\bfseries\upshape, fontupper=\slshape,
colback=red!10,colframe=red!75!black},
}

\numberwithin{equation}{section}
\usepackage[english]{babel}
\usepackage[bookmarks]{hyperref}

\usepackage{hyperref}
\hypersetup{
    colorlinks=true,
    linkcolor=blue,
    filecolor=magenta,      
    urlcolor=cyan,
    pdfpagemode=FullScreen,
}
\newcommand{\R}{\mathbb{R}}

\newcommand{\calplamb}{\mathcal{P}_{\Lambda}}

\date{}

\begin{document}
\title{Parameter Estimation for the Truncated KdV Model through a Direct Filter Method}
\author{
Hui Sun \thanks{Department of Mathematics, Florida State University}, 
\and Nick Moore \thanks{Department of Mathematics, Colgate University} ,
\and Feng Bao  \thanks{Department of Mathematics, Florida State University, \tt{bao@math.fsu.edu}} 
}
\maketitle

%\tableofcontents

\begin{abstract}
In this work, we develop a computational method that to provide real-time detection for water bottom topography based on observations on surface measurements, and we design an inverse problem to achieve this task. The forward model that we use to describe the feature of water surface is the truncated Korteweg–De Vries equation, and we formulate the inversion mechanism as an online parameter estimation problem, which is solved by a direct filter method. Numerical experiments are carried out to show that our method can effectively detect abrupt changes of water depth.
\end{abstract}

\textbf{keywords:} KdV equation, inverse problem, parameter estimation, data assimilation, particle filtering

\section{Introduction}

Extremely large surface waves, known as rogue waves, have been observed in oceans around the globe and studied from both theoretical and experimental perspectives. Of particular interest is when surface wave statistics deviate from Gaussian, or become anomalous. Studies have pointed to a variety of mechanisms that can produce anomalous wave statistics, one being an abrupt change in bottom topography. Through controlled laboratory experiments, Bolles et al.~\cite{BSM} demonstrated that strongly non-Gaussian surface-wave statistics can arise a short distance downstream of an abrupt depth change. These observations were explained by a theoretical model developed by Majda et al.~\cite{MMQ} based on statistical and dynamical analysis of the variable-coefficient truncated Korteweg–De Vries (KdV) system. This model accurately predicted a range of features observed in the experiments, including the transition in surface-displacement skewness as waves cross the abrupt depth change.

Given our understanding of how abrupt depth changes could influence the water surface waves and generate anomalous waves, in this work we consider the following inverse problem: How can we use observations on surface waves to infer the bottom topography such like the changes of depth? Such a problem may be of practical importance for technicians who seek to infer characteristics of bottom topography when observations are limited to surface wave measurements. More generally, this line of questioning applies to any situation in which one seeks to infer, or conceal, structure from dispersive wave scattering. Applications range from optical cloaking \cite{cai2007optical}, to detection of submarines from surface waves \cite{xue2020wake}, to the inference of black holes from gravitational waves \cite{gais2022inferring}.

\vspace{0.5em}
The bottom topography detection problem that we shall construct is based on the theoretical studies in \cite{MMQ, Moore_rewrite}, and we aim to detect and instantly infer the abrupt depth changes underwater based on the real-time measurements of surface waves. The mathematical tools that we adopt to solve this inverse problem include a so-called RK4 method for the tKdV model (introduced in \cite{MMQ}) and the optimal filtering techniques.

With the physics described by the tKdV model and the forward simulation implemented by the RK4 method, the key to solve the inverse problem is the development of an online estimation method for the parameters that reflect the abrupt changes in the tKdV model. To achieve this goal, we introduce a ``direct filer'' method to dynamically estimate the unknown parameters as we receive the observational data in the online manner. The direct filter method adopts the general framework of optimal filtering for hidden stochastic dynamical systems. The standard approach for the optimal filtering problem is the sequential data assimilation, which carries out Bayesian inference recursively to utilize the observational information to estimate the state of the hidden dynamical system. 

The main idea of the direct filter method is to combine the physics model with the observations on the state through the likelihood and use Bayesian inference as a projection tool to ``directly'' map the observational information to the parameter space. Our preliminary research shows that the direct filter is a very accurate online parameter estimation method for high dimensional data assimilation \cite{BSDE_KL_22, Bao_Direct}, and it's applicable to solve practical problems \cite{Bao_Bio, Bao_Atomic_2021}.

%However, one would see that estimating the depth change boils down to parameter estimation for the signal processes together with the coefficients in the tKdV which is typically of high dimensions (10 or more). 

%For this reason, the conventional filtering techniques like the augmented ensemble Kalman filters would not present good performance in general as pointed out in \cite{GRL}. To overcome this difficulty, we adopt the technique called the ``Direct Filter Method" proposed in \cite{GRL}.  As a highlight, we design a case study where the water wave encounters multiple step changes to demonstrate the efficiency of our method. 

The rest of this paper is organized as follows: In section 2, we introduce the physics background of the problem. In particular, we introduce our truncated KdV model, and we demonstrate that a standard RK4 method will generate stable state process with very little loss in the Hamiltonian. In section 3, we shall discuss the optimal filtering as a tool for online parameter estimation, and we will introduce our direct filter method with a particle implementation.
In section 4, we introduce how to apply the to estimate the unknown parameters in the rKdV model, and we shall carry out some numerical experiments to demonstrate the baseline performance of the direct filter in parameter estimation regarding the tKdV model. In Section 5, as a highlight, we design several experiments where abrupt depth changes in the bottom topography take place, and we will show that through parameter estimation we can dynamically detect the the hidden changes underwater.

\section{Preliminaries}

In this section, we present the physics background of our work by introducing the mathematical model, i.e. tKdV, motivated by the experiments conducted in \cite{BSM}. %Then we shall give a brief discussion on the development of the tKdV (truncated KdV) model. 

\subsection{KdV with variable coefficients}

\begin{figure}[hbt!]
  \centering
  \includegraphics[width=\linewidth]{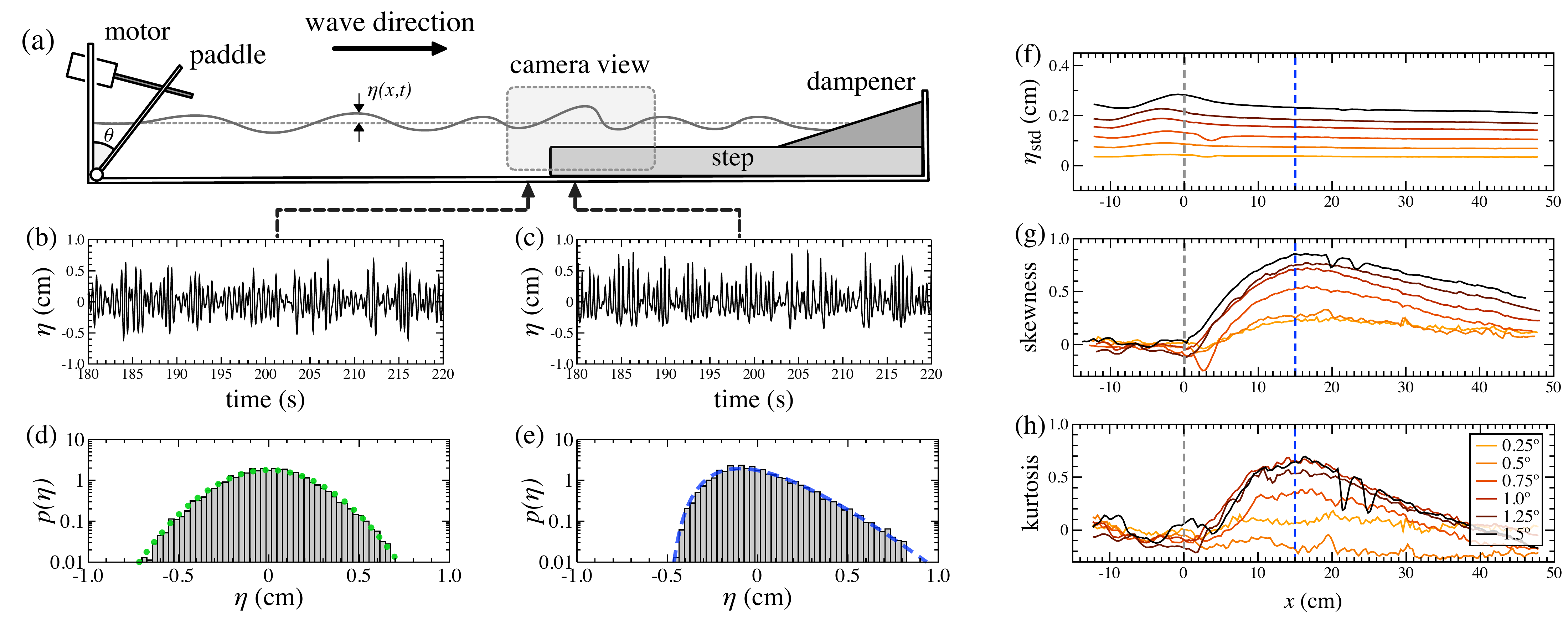}
  \caption{Laboratory experiments of Bolles et al.~\cite{BSM}. (a) Waves generated by a motor-paddle assembly propagate unidirectional through a long, narrow wave tank. Midway through the tank, the waves encounter an abrupt change in bottom depth which fundamentally alters their statistical distributions. Wave reflections are mitigated by a dampener placed at the end of the tank. (b)--(c) Surface-displacement measurements shown as a time series upstream (d) and downstream (c) of the abrupt depth change. (d)--(e) Histograms of the surface displacement measurements reveal Gaussian statistics upstream (b) and skewed statistics downstream (e) of the abrupt depth change. (f)--(h) Standard deviation, skewness, and kurtosis respectively of the surface displacement plotted against distance $x$ from the abrupt depth change. The depth change triggers a drastic rise in skewness and kurtosis, indicitave of strongly non-Gaussian statistics.
  }
  \label{fig:experiment}
  \end{figure}

The truncated KdV model (tKdV) introduced by Majda ~\cite{MMQ, Moore_rewrite} was developed to explain the experimental observations of Bolles et al.~\cite{BSM}. In those experiments, a field of randomized, unidirectional surface waves propagate through a narrow wave tank as illustrated in Fig.~\ref{fig:experiment} (a). Midway through the tank, the waves encounter an abrupt change in the bottom depth, which fundamentally alters their statistical distribution as the waves continue downstream. At the far end of the tank, reflections are mitigated by a dampener so that the wave-field remains nearly unidirectional. 

Figures~\ref{fig:experiment}(b)-(c) show time-series measurements of the surface displacement, $\eta$, a short distance upstream (b) and downstream (c) of the abrupt depth change. Figures~\ref{fig:experiment}(d)-(e) show the corresponding histograms of these measurements. A short distance upstream of the depth change, waves adhere to nearly Gaussian statistics as seen in Fig.~\ref{fig:experiment}(d). After passing over the depth change, however, the statistics skew considerable towards large positive displacement as seen in Fig.~\ref{fig:experiment}(e). This distribution contains a larger number of extreme events. In particular, Bolles et al.\cite{BSM} showed that the probability of a rogue wave  increases by a factor of roughly 50 downstream of the depth change.

Figures~\ref{fig:experiment}(f)-(h) show measurements of the first several moments of surfac displacement, including the standard deviation, the skewness and kurtosis. The measurements are displayed as a function of longitudinal position $x$, with $x=0$ corresponding to the location of the depth change; Different colors represent experiments with different driving amplitudes. The standard deviation of surface displacement provides the simplest estimate for the characteristic amplitude of waves. Interestingly, the standard deviation is relatively insensitive to longitudinal position, implying that the depth change does not substantially alter the typical size of waves. The skewness and kurtosis, however, change drastically as waves cross the depth change. Upstream, both quantities are nearly zero, consistent with nearly Gaussian statistics, but downstream of $x=0$ both skewness and kurtosis rise to a value on the order of unity, indicating a strong deviation from the Gaussian state.

To explain these experimental observations, Majda et al.~(2019) \cite{MMQ} proposed a shallow-water wave model based on the variable-coefficient, truncated KdV equation (tKdV). In this subsection, we provide a brief introduction to this theoretical framework, and we direct the interested readers to \cite{AM}, \cite{JB}, \cite{MMQ} and \cite{Moore_rewrite} for further information. 

The variable-coefficients KdV equation takes the form: 
\begin{equation} {\label{kdv_new}}
u_t + C_2 uu_x +C_3 u_{xxx}=0, \textit{ } x \in [-\pi, \pi]
\end{equation}
where the variable coefficients, $C_2$ and $C_3$, represent depth variations. We assume that such a system is normalized, with the mean and energy given by 
\begin{align}
	M&=\int_{-\pi} ^{\pi} u d x= 0 \\
E&=\frac{1}{2}\int_{-\pi} ^{\pi} u^2 d x= 1
\end{align}
This system possess the Hamiltonian:
\begin{equation}
H = \int^{\pi }_{-\pi } \frac{C_3}{6} u^3 - \frac{C_2}{2} \left( \frac{\partial u}{\partial x} \right)^2 dx
\end{equation}

Consider the situation of Bolles et al.~(2019) \cite{BSM}, in which waves pass from an upstream region $(x<0)$ with depth $d=d_-$ to a downstream region $(x>0)$ with depth $d=d_+$. We introduce a reference frame that travels with the leading-order wave speed, and choose the initial time so that the location of the depth change corresponds to $t=0$.
We introduce the dimensionless depth ratio
\begin{equation}
D=	\begin{cases}
		1 \ \ \text{for } \ t <0 \\
		D_+=d_+/d_- \label{crel}
	\end{cases}
\end{equation}
Then variable coefficients in KdV are given by
\cite{Moore_rewrite}:
\begin{equation}{\label{bigC}}
	C_2 : =c_2 D^{\frac{1}{2}}, \ \ C_3 := c_3 D^{-\frac{3}{2}}.
\end{equation} 
The coefficients $c_2$ and $c_3$ depend on experimental parameters, such as the frequency and amplitude of the wave forcing, and the leading-order wave speed. Importantly, these coefficients do not depend on the variable depth, and so, for any particular experiment, they remain constant in $x$ and $t$ \cite{Moore_rewrite}.

To numerically study the solution of this PDE, a Galerkin truncation is performed, i.e. the solution is written in the following form: 
$$ u(x,t) \sim \sum_{ |k| \leq \Lambda } \hat{u}_{k} e^{ikx} $$
Then the truncated system takes the following form:
\begin{equation}\label{tkdv}
\frac{d \hat{u}_{k}}{dt} +C_2 \mathcal{P}_{\Lambda}(uu_x)_k+C_3 \mathcal{P}_{\Lambda}(u_{xxx})_k=0
\end{equation} 
where $\calplamb$ is defined to be the projection with truncation number $\Lambda$. Numerical experiments show that a $\Lambda$ number of 16 will suffice \cite{nickhui}, \cite{MMQ}. And in later implementations, $\Lambda$ will be set equal to 16 throughout. This PDE can now be arranged to the following system of ODEs: 
\begin{equation}\label{rk4}
	\frac{d \hat{u}_k}{dt}=-C_3\frac{i k}{2}(\sum_{|k-m|\leq \Lambda} \hat{u}_{k-m} \hat{u}_{m}	) + i C_2 k^3\hat{u}_{k}.
\end{equation}
It is realized that the truncated system ($\ref{tkdv}$) also possesses the three conserved quantities (\cite{JB}, \cite{MMQ}, \cite{Moore_rewrite}) discussed in the previous section: 
\begin{enumerate}
	\item The Momentum: i.e. the zeroth mode of the spectrum $\hat{u}_0$ is 0.
	\item The Energy: $E = \pi \sum_{|k| \leq \Lambda} |\hat{u}_k|^2 = 2 \pi \sum_{0 < k \leq \Lambda} |\hat{u}_k|^2$
 	\item The Hamiltonian: $\mathcal{H}_{\Lambda}=C_3 H_{3,\Lambda}-C_2 H_{2,\Lambda}$
 \begin{equation}
 \textit{ }
H_{3,\Lambda} = \frac{\pi}{3} \sum_{m+n+l=0, |m|,|n|,|l| \leq \Lambda} \hat{u}_m\hat{u}_n\hat{u}_l, 
\textit{ }
H_{2,\Lambda} = \pi \sum_{|k| \leq \Lambda} k^2 |\hat{u}_k|^2
\end{equation}
\end{enumerate}
\subsection{Simulation strategies}
With the model given by the system of equations \eqref{rk4}, such equations need to be solved based on initial data and consequently, a time series of solutions will be obtained. Two techniques are in consideration. The first method is to use the symplectic integrator, which is adopted in \cite{MMQ}. Such method is designed to conserve both the energy and the Hamiltonian in long time horizon. However, one downside of this method is that it is computationally expensive since it is an implicit scheme. What is more, it will be more difficult to apply filtering method to such systems.    

\begin{figure}[hbt!]
  \centering
  \includegraphics[width=0.6\linewidth]{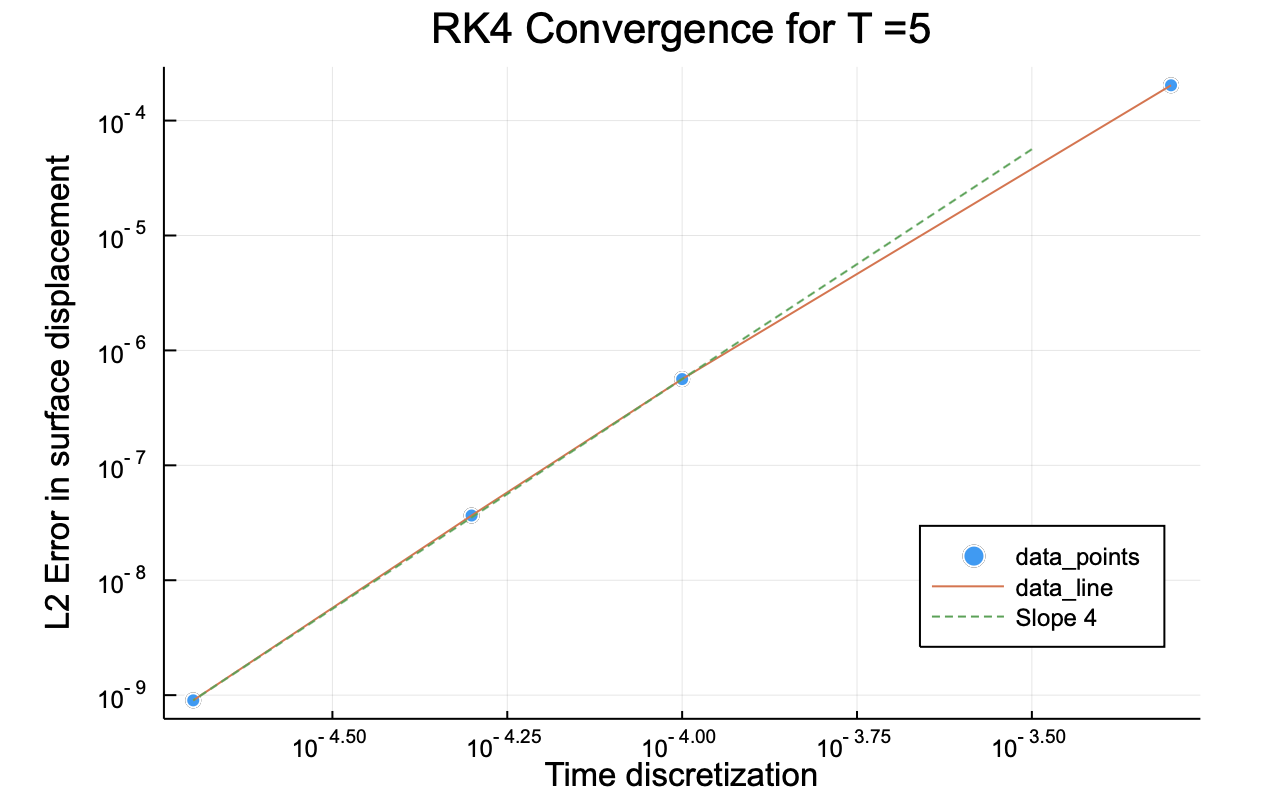}
  \caption{The log-log plot of error with respect to the time discretization, the straight yellow line has slope 4.0 }
  \label{fig:rk4_loss}
  \end{figure}

\begin{figure}[hbt!]
  \centering
  \includegraphics[width=0.6\linewidth]{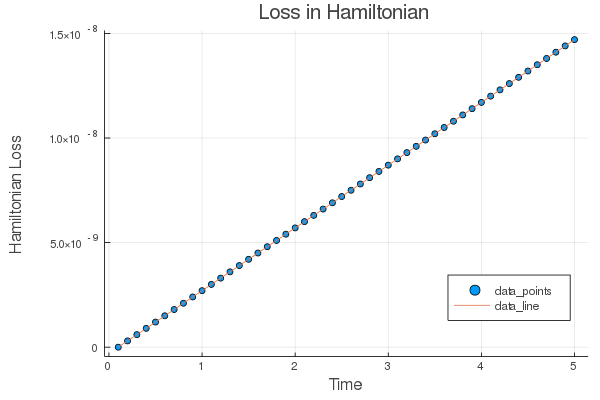}
  \caption{Loss in Hamiltonian}
  \label{fig:ham_loss}
  \end{figure}

Yet another idea will be to simply use RK4, which is a standard fourth order numerical scheme known to be stable over long time. And since the time horizon is not usually very large, the loss in Hamiltonian does not in general raise a concern. To justify the claim that RK4 will suffice, we plot figure \ref{fig:rk4_loss} and figure \ref{fig:ham_loss}.  By fixing initial data that follows the multidimensional Gaussian distribution, we let it evolve under the dynamics \eqref{rk4} with the RK4 as the numerical scheme. For all the systems tested below, we set the time horizon to be $T=5.0$. And the coefficients of the system are picked to be $C_2=1.0, C_3=1.0$.
\iffalse
\begin{figure}[!htb]
\center
  \begin{subfigure}[b]{0.5\textwidth}
    \includegraphics[width=\textwidth]{RK4_convergence.png}
    \caption{The log-log plot of error with respect to the time discretization, the straight yellow line has slope 4.0}
    \label{fig:2D21}
  \end{subfigure}
  %
  \begin{subfigure}[b]{0.5\textwidth}
    \includegraphics[width=\textwidth]{ham_loss.png}
    \caption{Error decay in iteration numbers}
    \label{fig:2D22}
  \end{subfigure}
  \caption{Loss in Hamiltonian}
  \label{2D_eg}
\end{figure}
\fi

By taking the temporal discretization level $\Delta t$ to be $10^{-6}$ as the benchmark result, i.e. the exact solution, we find that the log error between the numerical solution under different time discretization level and the exact solutions is roughly of order 4 in figure \ref{fig:rk4_loss}.
  
To see that long time behavior of the system \eqref{rk4} is stable and there is very minimum loss in the Hamiltonian of this system, we compute the Hamiltonian loss exactly. Given a fixed initial condition, one can compute the instantaneous Hamiltonian at time 0, and such value is used as benchmark since it should be conserved over time.  Then, to simulate \eqref{rk4}, we take $\Delta t=10^{-4}$, and for every $\delta t =0.1$, we compute the instantaneous loss in Hamiltonian to obtain figure \ref{fig:ham_loss}. 
  One can see that even though the loss in Hamiltonian gradually increases over time,  at the terminal time, the total loss in Hamiltonian is on the scale of $10^{-8}$ which is ignorable. And since the time horizon we take is never much larger than 5, it justifies our choice of choosing RK4 as the forward simulation method. 

\subsection{The Heatmap}
We observe from the experiments in Section 2.1 that the sudden change of the bottom topography of the tank will lead to more extreme water surface displacement statistics. Hence, in order to justify using \eqref{rk4} to model water propagation in a tank, we would do a direct simulation of the model and check that we can observe such phenomenon using the simulated data. 

\begin{figure}[hbt!]
  \centering
  \includegraphics[width=0.8 \linewidth]{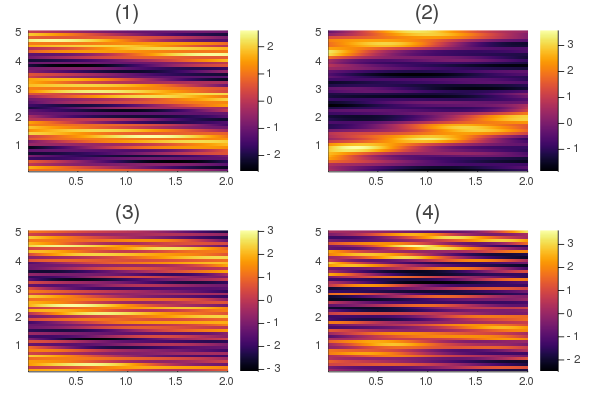}
  \caption{ (1) has $C_2=4.0, C_3=0.0$ and it is near Gaussian
 (2) $C_2=2.0, C_3=8.0$
 (3) $C_2=3.0, C_3=1.0$
 (4) $C_2=3.0, C_3=18.0$.}
  \label{fig:heatmap}
  \end{figure}

However, instead of using the standard statistical method (plotting histograms etc.) to achieve this goal, we would use direct simulation of the system and attempt to find anomalous waves triggered by the depth change. Now, we test four different systems with different level of nonlinearities by varying the coefficients $C_2, C_3$ (which in this case has nothing to do with $c_2,c_3$). We first test an extreme case where the nonlinearity term is absent by taking $C_2 = 4.0$ and $C_3=0$ whose plot is presented in Figure \ref{fig:heatmap} (a). In such heatmaps, the $X$ coordinate represents the space, the $Y$ coordinate represents the time and the color in the plot stands for the surface displacement.

 It is observed that the wave travels to the left in (a). And it makes sense due to the reasons as follows.
Since without the term involving $C_3$:
\[
u_t + C_2 u_{xxx}=0
\]
where $u_{xxx}$ is commonly known as the dispersive term.
By letting $u(x,t)=e^{i(yx-\sigma t)}$, one obtains the relationship: 
\begin{equation}{\label{dispersion}}
	-\sigma i -C_2 y^3 i=0 \Rightarrow \sigma=-C_2 y^3
\end{equation}
so by plugging \eqref{dispersion} back to $u(x,t)$, we obtain 
\begin{equation}
	u(x,t)=\exp \big(iy (x+C_2 y^2 t) \big)
\end{equation}
And so if the coefficient $C_2$ is positive, the wave should travels to the left which is exactly what we observed for this case

For the next stage, in Figure \ref{fig:heatmap} (3) we increase the value of $C_3$ to add more nonlinearity to the system by taking $C_2=3.0$ and $C_3=1.0$. However, the linear term is still dominant in this case, and so one should observe waves traveling to the left. 

Then in Figure \ref{fig:heatmap} (2)  we again increase both the magnitude of the $C_3$ term and the ratio of $C_3/C_2$ so that the nonlinear effect of the system becomes more pronounced. In this case, we pick $C_2=2.0$ and $C_3=8.0$, and the nonlinearity term starts to dominate, and different behavior from the previous two cases are observed: the wave starts to travel to the right. 

Finally in Figure \ref{fig:heatmap} (4), we make the nonlinearity effect more pronounced and we start to see more chaotic behaviors. 

This heat map is a collection of four connected systems, each solved individually using RK4. The initial data of the first system is given. And we use the output of the first system (The Fourier coefficients) at terminal time as the input of the second system and so on.

  \eqref{fig:heatmap}. The term $C_3$ indicates the nonlinearity in the fluids equation. And when this term grows larger, one will also observe waves that travels to the right in the heat map as shown in $(3),(4)$ in Figure \eqref{fig:heatmap}.

\begin{figure}
 \centering
  \begin{subfigure}[b]{0.8\textwidth}
    \includegraphics[width=\textwidth]{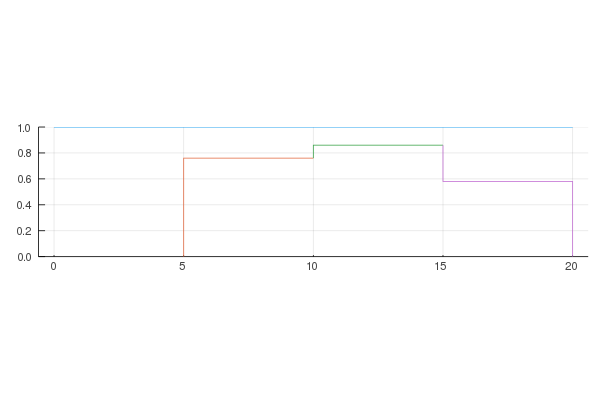}
    \vspace{-2.5cm}
    \caption{The distance between the free surface and the steps are 1.0, 0.24, 0.14,0.42}
    \label{fig:stepsplot}
  \end{subfigure}
  %
 % \vspace{4em}
  \begin{subfigure}[b]{0.8\textwidth}
    \includegraphics[width=\textwidth]{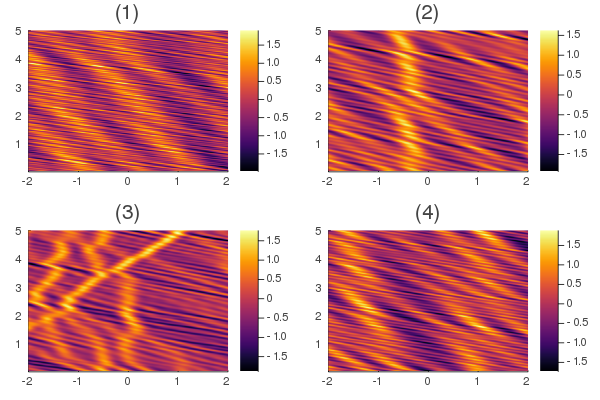}
    \caption{The plot of the surface displacement. It shows the surface change in space according to the time evolution}
      \label{fig:stepheat}
  \end{subfigure}
  \caption{}
\end{figure}

In order to further understand the behavior of wave statistics after passing over abrupt depth changes, we design a water tank with the bottom topography generally depicted in Figure \ref{fig:stepsplot}. The blue straight line at the top represents the free surface and the orange line stands for the first step. The step ratio is given by $d_+/d_-=0.24$ where we pick $d_-$ to be the $1.0$ which is the free water surface, and $d_+$ to be $$1- \frac{\text{height of step } }{\text{depth of free surface}}$$
In stead of letting the water wave pass over the step change once, we add two more steps following the first one and they are represented by the green line and the purple lines in  Figure \ref{fig:stepsplot} whose step ratio $d_+ /d_-$ are given by 0.14 and $0.42$ accordingly. 

To see how the sequential change of steps affect the wave statistics, we fix an initial state $\hat{u}_0$, let it evolve under the dynamics \eqref{rk4} with $C_2,C_3$ defined by \eqref{bigC} with $d_+/d_-=1.0$, and set the total time to $T_1=5$. After the inverse Fourier transform, we obtain the wave surface displacement and plot part (1) in Figure \ref{fig:stepheat}. One can observe that since in this case there is no step present, the $C_2$ term is more dominant and the wave is roughly linear. 

 Then, at the terminal time $T_1=5.0$, we take the terminal state of $\hat{u}_{T_1}$ and make it the initial state for the second step. Now, due to the fact that $d_+/d_-=0.24$, we have a new set of coefficients $C_2, C_3$ by \eqref{bigC}, and due to the fact that now the water is much shallower than the previous case, in part (2) of Figure \ref{fig:stepheat}, we observe that waves of large magnitude starts to show up and there are more chaotic behaviors in such system. 
 
 By repeating the same procedure, we also obtain Figure \ref{fig:heatmap} (3) and Figure \ref{fig:heatmap} (4). We point out that for Figure \ref{fig:heatmap} (3), we are in a scenario where the water is very shallow, and it is seen in the plot that there is a wave of large magnitude that travels to the right, which is some phenomenon that one does not observe in Figure \ref{fig:heatmap} (1) or Figure \ref{fig:heatmap} (4).

\begin{figure}[hbt!]
  \centering
  \includegraphics[width=0.6 \linewidth]{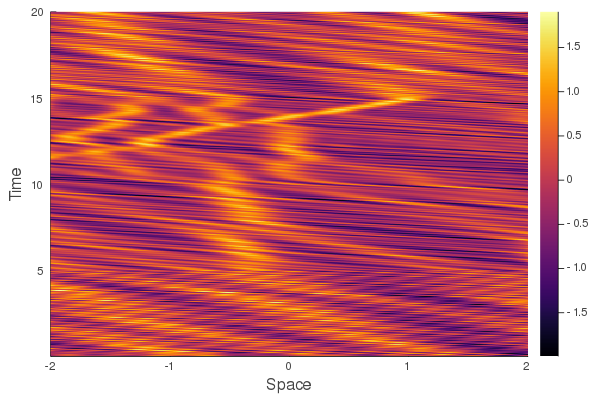}
  \label{fig:stackheatmaps}
   \caption{The stacked heatmaps of water waves in different region }
  \end{figure}

\vspace{2em}

The forward problem described above helps us understand the water surface displacement statistics, and we find that with an abrupt depth change from a deep water region to a shallow water region it is more likely for us to observe waves of much larger amplitude.  
Based on the physics knowledge we discussed above, in this work, we aim to study the inverse problem of the KdV equation, and we consider the scenario in which by observing the water surface displacements, which can be simulated by the tKdV model, how to find the depth of hidden abrupt changes underwater.  In what follows, we introduce a direct filter method for online parameter estimation to achieve this goal.

\section{A direct filter method for parameter estimation}
%The forward problem described above helps us understand the water surface displacement statistics, and we find that with an abrupt depth change from a deep water region to a shallow water region it is more likely for us to observe waves of much larger amplitude.  

In this section, we are interested in studying the following inverse problem: given the online (RK4) solution of the tKdV, what can we say about the bottom topography of the water tank (the level of the depth change, the location of the depth change etc.)? Such a problem is of practical interest, because one typically faces the following filtering problem: given realtime data, what kind of useful information can we extract from it? For example, in our case, given the incoming water surface displacement data, can we tell when will we have a sudden change of bottom topography?

To answer such questions, we need to adopt some techniques from data assimilation. 

\subsection{General setup of the optimal filtering problem}
There are two main problems in data assimilation, the first one is the smoothing problem which is typically an offline problem. One usually needs to have a very large dataset and the simulation process needs to be done all at once. For this method, it is required that all the data being collected for investigation, and so this is a ``static process".  The other method, called optimal filtering, is an online method. In an optimal filtering problem, one can pick a time window and data will be collected as time evolves. Since such method is more realistic for our purpose, we are going to study this approach closely. 

Suppose one has the following signal process $S_n \in \mathbb{R}^d$, 
\begin{equation}\label{non_linear_s}
	S_{n+1}=g(S_n) + w_n, \hspace{1em} n = 0, 1, 2, \cdots,
\end{equation}
where $g:\R^d \rightarrow \R^d$ is a given dynamical model, and the term $w_n \in \R^d$ is the noise in the system which is usually assumed to be Gaussian. In the case where $w_n$ is absent and $g$ contains a $\Delta t$ term, such system of equations is deterministic can be treated as ODEs. 

However, the real signal process $S_n$ is often not directly observed, and the observational data one holds are typically the transformation of the real process: 
\begin{equation}
	O_{n+1}=G(S_{n+1})+\xi_{n+1} 
\end{equation}
where we use $O_n \in \R^k$ to denote the  observational data, and $G: \R^d \rightarrow \R^k$ could be either a linear or a nonlinear function, and $\xi_{n+1}$ is Gaussian noise. 
Notice that the observed process can have different dimensions than the real process. 

Now, the goal of optimal filtering is to find the conditional expectation:
$$ \mathbb{E}[S_{n+1}|O_{1:n+1}]$$
which has the following interpretation: given the observational data from time $1$ to $n+1$, what is the best estimate of the real system state $S_{n+1}$? 

In this work, we use the Bayesian inference framework to obtain such estimation. And the filtering process takes two steps. The first step is the ``prediction" which is realized by the Chapman-Kolmogrov formula
\begin{equation}{\label{predict}}
	p(S_{n+1}| O_{1:n}) = \int p(S_n+1{}|S_n) p(S_n|O_{1:{n}})dS_n
\end{equation}
where we assume that $p(S_n|O_{1:{n}})$ is known, and $p(S_n+1{}|S_n)$ is the transition kernel. 

The next step uses the Bayes' formula, and will give us the update of the new posterior:
\begin{align}
	p(S_{n+1}|O_{1:n+1}) &=\frac{p(O_{n+1}|S_{n+1},O_{1:n})p( S_{n+1}|O_{1:n})}{p(O_{n+1}| O_{1:n})} \\ 
	&=\frac{p(O_{n+1}|S_{n+1})p( S_{n+1}|O_{1:n})}{p(O_{n+1}| O_{1:n})}
\end{align}

\begin{remark}
Here $p(O_{n+1}|S_{n+1})$ is the likelihood function. We dropped the dependence on $O_{1:n}$ since we already have the information for $S_{n+1}$ suffices due the Markov property of the system. $p( S_{n+1}|O_{1:n})$ comes from the prediction step and the denominator is the normalizing constant of the numerator. 
\end{remark}

For the linear model, Kalman filters are usually used for the above steps. However, for nonlinear problems, the particle methods are adopted since not all the analytic expressions for each term are available. The general idea of such methodology is that one constructs an ensemble of particles (a convex combination of Dirac measures) and update their locations and weights accordingly so that in the limit of large particle numbers, the conditional distribution will converge to the exact one. 

Let's call those particles $\lbrace s_n^{m} \rbrace_{m=1}^M$. 
By using the empirical distributions, we have the following results:
\begin{enumerate}
	\item $\tilde{p}(S_{n}|O_{1:n}) :=\sum_{m=1}^M w^{(m)}_{n}\delta_{s_n}^{(m)} (S_n)$ where $w^{(m)}_n$ stands for the weight of the particle $m$. 
	\item The prediction step consists of obtaining the following empirical distribution: $\tilde{\pi}(S_{n+1}|O_{1:n})=\sum_{m=1}^M w^{(m)}_{n} \delta_{\tilde{s}_{n+1}}^{(m)} (S_{n+1})$ which is an approximation for $p(S_{n+1}| O_{1:n})$, and it is obtained from the Chapma-Kolmogrov formula. $\tilde{s}_{n+1}$ stands for the empirical particle updated after one step, and they are obtained by following the dynamics \eqref{non_linear_s}.  Notice that now we have changed the particle positions but not their individual weight.  
	\item The empirical distribution of the update step is given by the following
	\begin{equation}{\label{e_update}}
		\tilde{\pi}(S_{n+1}|O_{1:n+1}):= \frac{\sum_{m=1}^M p(O_{n+1}|\tilde{s}^{(m)}_{n+1}) w^{(m)}_{n}\delta_{\tilde{s}_{n+1}}^{(m)} (S_{n+1})}{\sum_{m=1}^M w^{(m)}_{n} p(O_{n+1}|\tilde{s}^{(m)}_{n+1})}
	\end{equation}
	
The denominator of the above equation is a constant and so we rewrite \eqref{e_update} in the following manner: 
\begin{align}
	\tilde{\pi}(S_{n+1}|O_{1:n+1})=\sum_{m=1}^M w^{(m)}_{n+1}\delta_{s_{n+1}}^{(m)}\\  
	\tilde{w}^{(m)}_{n+1}=p(O_{n+1}|\tilde{s}^{(m)}_{n+1}) w^{(m)}_{n}\\
	w^{(m)}_{n+1}=\tilde{w}^{(m)}_{n+1}/\sum^M_{n=1}\tilde{w}^{(m)}_{n+1}
\end{align}

	\item To avoid degeneracy (meaning that a lot of the weights of the particles are ignorable), one needs a resampling step for the updated measure. And we do this by generating M samples from the distribution (\ref{e_update}) and assign weight $\frac{1}{M}$ to each of the weight. 
	Again, the goal of the step is to kill the small weights. (One does not want to avoid sampling from the tail of the distribution.) 
\end{enumerate}
We point out that it is important to use particle filters in our problem setup since there is no Gaussian assumptions made on our system. However, it is widely recognized that in order to avoid the degeneracy issue, an exponential growing number of particles are needed as the dimension increases. In our case, the complex system of interest \eqref{rk4} is of size $\Lambda$, and so there will be $2 \Lambda + 2$ unknowns. Since typically $\Lambda$ will be 10 or larger, the dimension for our setup will surely lead to degeneracy. 

On the other hand, it is observed that only the two parameters $(C_2, C_3)$ in \eqref{rk4} are of major interest and the state estimates are simply unnecessary. As a result, we will adopt the \textit{Direct Particle Filter} method proposed in \cite{Bao_Direct}. Such method will gives us the benefits of: 
\begin{enumerate}
	\item The improvement in algorithm efficiencies, since the dimension of the problem has been significantly reduced, and so the degeneracy problem can be overcome. 
	\item The procedure of estimating the state process has been avoided, and so the algorithm will only focus on the two parameters of interest
\end{enumerate}
In the next subsection, we will describe the setup and implementation procedure of the direct filter.

\subsection{Parameter estimation: the direct filtering method}
Since we aim to estimate the parameter values in the tKdV model, we will use $\theta$ to denote the parameter of interest and write down the following system for the parameter estimation problem: 
\begin{align}
	X_{n+1}&=h(X_n,\theta)+w_n \label{state} \\
	Y_{n+1}&=HX_{n+1}+\xi_{n+1} \label{estimate}
\end{align}
where $h: \R^d \times \R^p  \rightarrow \R^d$ is a nonlinear function representing the physics model that we are considering; $H: \R^d \rightarrow \R^m$ is a linear matrix; $(w_n, \xi_n )$ are independent Gaussian noises. In the above dynamical system, $X$ describes the state of some physics model, and $Y$ provides direct observations on state $X$ with noise perturbation $\xi$. The parameter estimation problem and we are interested in is to estimate $\theta$ in \eqref{state} by using the observational data $Y$ in \eqref{estimate}.  
%\iffalse
%The augmented filter strategy of dealing with the parameter estimation is to treat the two objects $\theta$ (the parameter to be estimated) and $X_{n}$ (the real process) as a pair and use sequential Bayes' inference to estimate $\theta$ under the state space model.  That's, formulating the following processes:
%$$F(S_n)=(\theta_n, h(X_n,\theta_n))$$
%and study the following joint posterior distribution $p(S_{n+1}|O_{1:n+1})=p(\theta_{n+1},X_{n+1}|O_{1:n+1})$. 
%\fi
In what follows, we give a brief discuss for the direct filtering method. The more rigorous formulation and discussion can be found in \cite{Bao_Direct}.

In the direct filter, instead of treating the parameter $\theta$ as a deterministic constant, we consider $\theta$ as a stochastic process to be estimated with respect to time.  To proceed, $\theta$ in \eqref{state} is replaced with $\theta_{n}$, and we rewrite \eqref{state} and \eqref{estimate} in the following form:
\begin{align}
	\theta_{n+1}&= \theta_n + \epsilon_n\\
	Y_{n+1}&= H (h(X_n, \theta_{n+1})+w_n)+\xi_{n+1} \label{new-strat}
\end{align}
We define $\eta_{n+1}=H w_n +\xi_{n+1}$, which is a is Multivariate gaussian variable, and we can obtain the following dynamics:
\begin{align}
	\theta_{n+1}&= \theta_n + \epsilon_n \label{new-strat1} \\
	Y_{n+1}&= H h(X_n,\theta_{n+1})+\eta_{n+1}. \label{new-strat2}
\end{align} 
The major difference between the above formulation and the standard formulation \eqref{state}-\eqref{estimate} of the parameter estimation problem is that instead of treating $\theta$ as a constant, we introduce a pseudo dynamics to formulate $\theta$ as a stochastic process. As a result, we aim to find $\mathbb{E}[\theta_n | Y_{1:n}]$ as our ``best'' estimate for $\theta_n$ at any time instant $n$, where $Y_{1:n}$ contains the information of the observation process $\{Y_i\}_{i=1}^n$.

%If one only considers $\theta$ as the target state process in the optimal filtering problem, the auxiliary state process $X_n$ in \eqref{new-strat2} should be either eliminated or approximated by a function of $(Y_n, \theta_n)$. And we make the approximations by the following considerations. 

To incorporate the observational data $Y_{n+1}$ to estimate $\theta_{n+1}$, we apply the Bayes formula to obtain the posterior distribution as follows
\begin{equation}\label{Bayes}
p(\theta_{n+1}| Y_{1:n+1}) =  \frac{p(Y_{n+1} | \theta_{n+1}) p(\theta_{n+1} | Y_{1:n})}{C},
\end{equation}
where $C$ is a normalization factor (i.e. the marginal likelihood), $p(\theta_{n+1} | Y_{1:n})$ is the prior distribution for the target parameter at time $n+1$. Since the state dynamics for the target parameter variable is a zero dynamics, the derivation of the prior distribution is straightforward, and the main effort lies on evaluating the likelihood $p(Y_{n+1} | \theta_{n+1})$. 

Note that the original physics model $h$ in \eqref{state} is the bridge that connecting $\theta$ and $X$, where $X$ is observed through $Y$, we have the following expression for the likelihood function
\begin{align}
	p(Y_{n+1}|\theta_{n+1}) &= \int p(Y_{n+1}|X_{n+1})p(X_{n+1}|X_n, \theta_{n+1})dX_{n+1}\label{Bayes}%\\
	%&=\int (\int p(Y_{n+1}|X_{n+1})p(X_{n+1}| X_n, \theta_{n+1}) p(X_n|\theta_{n+1}) dX_n) dX_{n+1}
\end{align}
%Approximating $X_n= H^{-1}(Y_n -\xi_n) \approx H^{-1} Y_n$, we have the following: 
%Note that we don't have direct information about $X_n$. Instead, we have the observational data $Y_n$ on $X_n$, and we can use the following approximation 
Different from other online parameter estimation methods, instead of generating a long-term simulation trajectory for the process of $X_n$ (e.g. \cite{KDSMC}), we use the fact that the observational data $Y$ provides direct observations on $X$ and introduce the following approximation scheme
$$H^{-1} Y_n \approx H^{-1}(Y_n -\xi_n) = X_n$$ 
to approximate $X_n$. Then, we can estimate the likelihood as
\begin{equation*}
	p(Y_{n+1}|\theta_{n+1}) \approx \int p(Y_{n+1}|X_{n+1})p(X_{n+1}| H^{-1}Y_n, \theta_{n+1})  dX_{n+1}.
\end{equation*}	
%The following form shows the strategy of approximation
%\begin{align}
%	\theta_{n+1}&= \theta_n + \epsilon_n\\
%	\tilde{Y}_{n+1}&= H h(H^{-1}\tilde{Y}_{n},\theta_{n+1})+\eta_{n+1} \label{new-strat}
%\end{align} 
%where $\tilde{Y}_n = Y_n$ which are the original observational datasets.
%
%Such formulation of the problem potentially avoids the curse of dimension problem since now the states of interest are just the parameters themselves. In the next section, we will outline the implementation strategies. 
%\begin{remark}
%	Refer to \cite{GRL} for some comments about the approximation on page 7 of the paper.  
%\end{remark}
In what follows, we introduce a particle implementation of the direct filter methodology. More detailed discussions can be found in \cite{Bao_Direct}

%\subsubsection{Particle implementation of the direct filter}

To initialize our algorithm, we first generate a collection of $M$ particles, denoted by $\lbrace \theta_0^{(m)}\rbrace_{m=1}^M$, from the initial guess for the parameter to be estimated. 

From time instant $n$ to $n+1$, with a set of particles $\lbrace \theta_n^{(m)}\rbrace_{m=1}^M$ that describe $p(\theta_{n}|Y_{1:n})$ at time instant $n$, we implement a ``prediction step'', an ``update step'', and a ``resampling step'' as what follows:

\begin{enumerate}
	\item \textbf{Prediction step.}  The prediction step generates a prior estimate for the target parameter $\theta_{n+1}$ before receiving the new data.
	
	Specifically, we add $\lbrace \epsilon_n^{(m)}\rbrace_{m=1}^M$ to $\lbrace \theta_n^{(m)}\rbrace_{m=1}^M$ through \eqref{new-strat1} to get a set of predicted particles, i.e. 
$$ \tilde{\theta}_{n+1}^{(m)} =  \theta_n^{(m)} + \epsilon_n^{(m)}, \hspace{1em} m = 1, 2, \cdots, M.$$ 

The prediction step gives the following empirical distribution for the prior $p(\theta_{n+1}| Y_{1:n})$ through the particles $\lbrace \tilde{\theta}_{n+1}^{(m)}\rbrace_{m=1}^M$ as 
	\begin{equation}\label{MC:prior}
		\tilde{\pi}(\theta_{n+1}| Y_{1:n})=\frac{1}{M}\sum_{m=1}^M \delta_{ \tilde{\theta}_{n+1}}^{(m)} (\theta_{n+1}).
	\end{equation}
	
	\item \textbf{Update}. The update step incorporates the observational data and derive a weighted posterior distribution based on the prior derived in \eqref{MC:prior}, i.e.
	\begin{equation}\label{MC:Bayes}
		\tilde{\pi}(\theta_{n+1}| Y_{1:n+1})=\frac{\sum_{m=1}^M p( Y_{n+1}|\tilde{\theta}^{(m)}_{n+1}) \delta_{ \tilde{\theta}_{n+1}}^{(m)} (\theta_{n+1})}
		{C},
	\end{equation}	
and the likelihood function in \eqref{MC:Bayes} is given by
\begin{equation}\label{og_likelihood}
	p(Y_{n+1}|\tilde{\theta}^{(m)}_{n+1}) = \exp( -\frac{1}{2} ||Hh(H^{-1}Y_n, \tilde{\theta}_{n+1}^{(m)})- Y_{n+1}||^2_R),
\end{equation}
where $||\alpha||_R:=\alpha R^{-1} \alpha$ with $R$ standing for the invariance variance of the observational noise $\eta_{n+1}$.

By choosing the weight for each particle $\tilde{\theta}^{(m)}_{n+1}$ as the normalized likelihood, i.e. $\omega^{(m)}_{n+1} : = p(Y_{n+1}|\tilde{\theta}^{(m)}_{n+1})/C$, we obtain a weighted empirical distribution characterized by $\{(\tilde{\theta}^{(m)}_{n+1}, \omega^{(m)}_{n+1})\}_{m=1}^M$.

\item \textbf{Resampling step.}
The purpose of the resampling step is to generate a set of equally weighted samples to avoid the degeneracy issue \cite{BSDE_filter, BSDE_filter_15, BaoCC_2019, BaoCH_CiCP, Bao_CMS}.

In this work, we simply use the importance sampling method to generate samples, denoted by $\lbrace \theta_{n+1}^{(m)}\rbrace_{m=1}^M$ , from the weighted importance distribution $\tilde{\pi}(\theta_{n+1}|Y_{1:n+1})$ and replace particles with small weights by the particle locations with large weights \cite{MDGPFB}. 
\end{enumerate}

With the ``prediction - update - resampling'' procedure, the estimate that we obtain for the target parameter at time instant $n+1$ is given by 
\begin{equation}
	\tilde{\theta}=\frac{1}{n+1-l} \sum_{i=l}^{n+1} \tilde{\mathbb{E}}[\theta_i | Y_{1:i}],
\end{equation}
where $\tilde{\mathbb{E}}[\theta_i | Y_{1:i}]$ is the mean estimate of the empirical distribution $\tilde{\pi}(\theta_{i}|Y_{1:i})$ obtained at the time instant $i$.
Note that $l$ is a number of burn-in steps to reduce the influence of large noises at some time instants.  

In the Appendix of this paper, we use an illustration example to demonstrate the baseline effectiveness of the algorithm. 
 
\section{The direct filter approach for parameter estimation for the tKdV model}\label{DF-tKdV}
In this section, we shall explain how to apply the direct filter to estimate the parameters $C_2$ and $C_3$ in the tKdV model, and we demonstrate, through numerical experiments, that $C_2$ is a more effective indicator for the depth of the shallow water. 

%Through our numerical experiments, we claim that it's preferable to estimate the parameter $C_2$ than $C_3$ due to the fact that $C_2$ is a more sensitive parameter over $C_3$.
 
\subsection{Problem setup}
The goal of parameter estimation in the tKdV model is to find the up stream/down stream depth ratio: 
\begin{align}
	D=d_+/d_- 
\end{align}
given the pre-determined constants $c_2=0.0236$ and $c_3=0.1965$ as studied in \cite{Moore_rewrite}. To this end, we propose to precisely estimate the coefficient $C_2=c_2D^{\frac{1}{2}}$ and/or $C_3=c_3D^{-\frac{3}{2}}$ both in the upstream and downstream region.

We assume that the signal process, i.e. $\{X_n\}_n$ in \eqref{state}, is described by the solution of equation \eqref{rk4}, and we choose $\Lambda=16$ in this model. Hence, denoting the RK4 operator by $\Phi$, we have 
\begin{equation} \label{rk4_data}
	\hat{u}_n^k=\Phi^k(\hat{u}^1_1, ..., \hat{u}_n^{\Lambda}),  \ \forall k=1,...,\Lambda,
\end{equation}
where $k$ denotes the $k-th$ component of the state vector $\Phi(\hat{u}^1_1, ..., \hat{u}_n^{\Lambda})$, and $n$ is the time index. In this way, we let the state process be the solution vector $\hat{u}_n\in \mathbb{C}^{\Lambda}$, i.e.  $X_n = \hat{u}_n$, and we choose $\Delta t=10^{-4}$. Hence, we will have a total collection of $\mathcal{T}=\frac{T}{\Delta t}$ such state vectors, where $T$ is the terminal time. 
%\begin{remark}
%	We point out here that one does not need to store the entire $M$ vectors, since typically the historical data far back will not impact the filtering process as much as the current data, and we will pick a window of length $L$. 
%\end{remark}

To formulate the data assimilation problem for parameter estimation, we assume that we receive datasets with linear dependence on $\hat{u}_n$: 
\begin{align}
	Y_{n+1}=H \hat{u}_{n+1}+\xi_{n+1}
\end{align}
as the observational process, and the added noise, $\xi_{n+1}$, follows a Gaussian distribution. We denote $\theta : =(C_2, C_3)$ as our parameter of interest, and we add noise to $\theta$ to make it a stochastic process. The introduction of additional noise to the parameter vector of interest will transform them into an ensemble. As such, filtering/data assimilation techniques such as particle methods can be applied to facilitate the state estimation of those parameters. Since the KdV model describes waves on shallow water surfaces, it's a reasonable assumption to obtain direct observations on the solution, i.e. the surface waves.
As a result, we have the following optimal filtering problem for the parameter estimation task: 
\begin{align}
	\theta_{n+1}&= \theta_n + \epsilon_n \label{parameq} \\
	Y_{n+1}&=H \Phi( H^{-1}(Y_n), \theta_{n+1} ) + \xi_{n+1}. \label{obeq}
\end{align}
%For simplicity, we assume that $H:=I$ which is the identity matrix and\
%\begin{equation} \label{obs_noise}
%	\xi_{n}\sim 0.01 \mathcal{N}(0, I)
%\end{equation} 
%As a result, \eqref{obeq} is simplified to 
%\begin{equation}\label{obeq_simp}
%	Y_{n+1}= \Phi(Y_n, \theta_{n+1} ) + \xi_{n+1}
%\end{equation}
%As a result, \eqref{og_likelihood} takes the following form:
%\begin{align}
%	p(\tilde{Y}_{n+1}|\tilde{\theta}^{(m)}_{n+1}) &= \exp( -\frac{1}{2} ||\Phi(\tilde{Y}_n, \tilde{\theta}_{n+1}^{(m)})-\tilde{Y}_{n+1}||^2_R)\\
%	&= \exp( -\frac{10^4}{2} ||\Phi(\tilde{Y}_n, \tilde{\theta}_{n+1}^{(m)})-\tilde{Y}_{n+1}||^2)
%\end{align}

In what follows, we carry out some numerical experiments to study the applicability of using the direct filter method to estimate $C_2$ and $C_3$, which is equivalent to estimating the water depth parameter $D$.

\subsection{Numerical experiments to study the performance of the direct filter in estimating $C_2$ and $C_3$}
Recall that if one can obtain a good estimate for either $C_2$ or $C_3$, then the step ratio $D$ can be computed through relation \eqref{bigC}, which is a major goal of parameter estimation for the tKdV model.  

%In this subsection, we will design various parameter setups and present their implementation results accordingly, we will compare the accuracy and stability results and comment on their behaviors.  We first carry out numerical experiments to estimate the value of $C_2$, and we shall discuss the parameter estimation performance for $C_3$ in a later experiment.

In our first numerical experiment in this section, we let  $\xi_{n}= \mathcal{N}(0,0.01^2)$ and choose $M=2000$ particles to empirically the distribution of the unknown parameter, and we introduce a total of $2500$ data assimilation steps. Moreover, we pick $\epsilon_n \sim \mathcal{N}(0,\Sigma)$, where $\Sigma =diag(0.3^2,0.017^2)$, and the observation matrix is chosen as $H = Id$.

\begin{figure}[hbt!]
  \centering
  \includegraphics[width=0.6\linewidth]{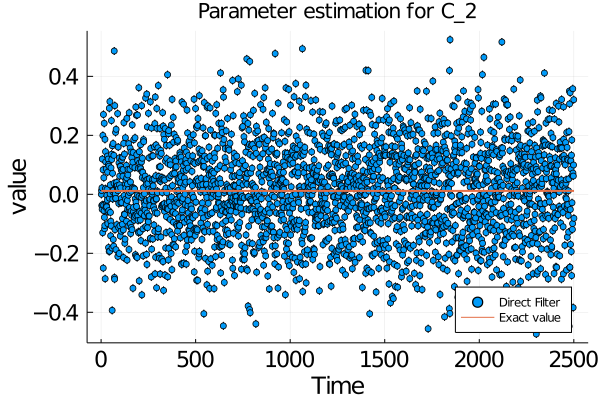}
  \caption{The estimation for paremeter $C_2$ }
  \label{fig:sd}
  \end{figure}

To generate the ``synthetic data'' for the state, we choose a pre-defined depth value $D = 0.24$, which will give us the real parameter, i.e. $C_2=0.01158$ and $C_3 = 1.671$. In Figure \ref{fig:sd}, we present the performance of parameter estimation by using the direct filter in estimating $C_2$, where the red line indicates the exact value for $C_2$, and the blue dots are the estimated value accordingly at each time. As one can observe from Figure \ref{fig:sd}, all the estimated parameter values are close to the real parameter value, and the variance of those estimates is $0.018$. Such an estimation variance is relatively large compared to the mean, but it's actually very small in terms of the absolute value.  By dropping $l = 400$ burn-in estimated values, we get the mean estimate for parameter $C_2$ as $\tilde{C}_2 = 0.0113$

\begin{figure}[hbt!]
  \centering
  \includegraphics[width=0.6\linewidth]{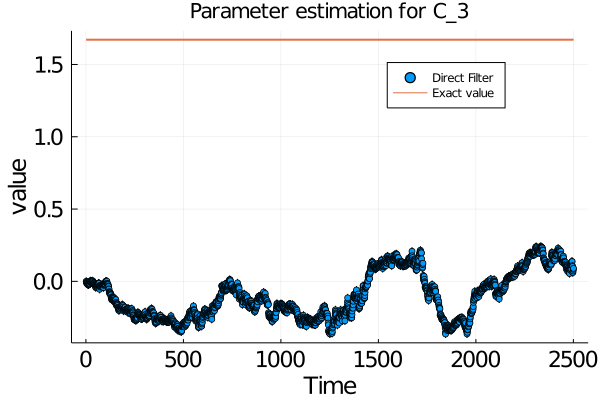}
  \caption{The estimation for paremeter $C_3$. }
  \label{fig:sd3}
  \end{figure}
  In Figure \ref{fig:sd3}, we present the parameter estimation results for $C_3$. From this figure, we can see that that the estimated values are not close to the true parameter value. In other words, the direct filter method can successfully estimate $C_2$, but it has difficulties in estimating $C_3$.

To further demonstrate the performance of the direct filter method in estimating $C_2$, we choose different values of $D$ and show the performance of parameter estimation for $C_2$, and the estimation results are presented in Table \ref{Values:D}, where $D = 1.0$ means the upstream region with no depth change took place; $D = 0.42$ and $D = 0.24$ are the values used in \cite{Moore_rewrite} throughout; and $D = 0.14$  corresponds to the case where the water is rather shallow. 
\begin{table}\caption{Performance of the direct filter with different choice of $D$}\label{Values:D}
\centering%\begin{center}
 \begin{tabular}{||c c c||} 
 \hline
 $D$ & $C_2$ & Estimated $C_2$  \\ [0.5ex] 
 \hline\hline
 1 & 0.236 & 0.234 \\ 
 \hline
0.42 & 0.0153  & 0.0155 \\
 \hline
 0.24 & 0.0115 & 0.0113 \\
 \hline
 0.14 & 0.0088 & 0.0095\\ [1ex] 
 \hline
\end{tabular}
%\end{center}
\end{table}
As one can see from the table, all the estimated values for $C_2$ are close to the true $C_2$ values in the mildly shallow or shallow regions ($D=1,0.42,0.24$).  And those values will in turn give us a good estimated value for $D$ (the depth change ratio) since $c_2$ is a known value. One may also observe that when we reach a region that is very shallow ($D=0.14$), the estimated $C_2$ will present a larger error. This could be explained by the fact that when $C_3$ is significantly larger than $C_2$, the system becomes rather chaotic, and it is hard for the optimal filter to learn the behavior of the linear term accurately. We also point out that it is rare sometimes to observe a very large depth change $D=0.14$ in reality. 

% \begin{figure}[hbt!]
%  \centering
%  \includegraphics[width=0.6\linewidth]{estimate_c3_024.png}
%  \caption{The estimation for paremeter $C_3$. It is observed that the predicted value is nowhere close to the true parameter value of interest. Or it may be that the burn-in time is too long for this particular parameter setup.}
%  \label{fig:sd3}
%  \end{figure}
 
 \vspace{0.5em}

%However, the estimation for $C_3$ presents some difficulties largely due to the fact that it requires much longer burin-in time and the oscillation over the mean is a lot more drastic than the case for $C_2$. Also, it is observed that if one wants to attain a reasonable prediction for $C_3$, the noise added $\epsilon_n \sim \mathcal{N}(0,\Sigma)$ should have different $\Sigma$ rather than taking $\Sigma=\text{diag}(0.3^2,0.017^2)$ as the following plot demonstrates. 
 \begin{figure}[hbt!]
  \centering
  \includegraphics[width=0.6\linewidth]{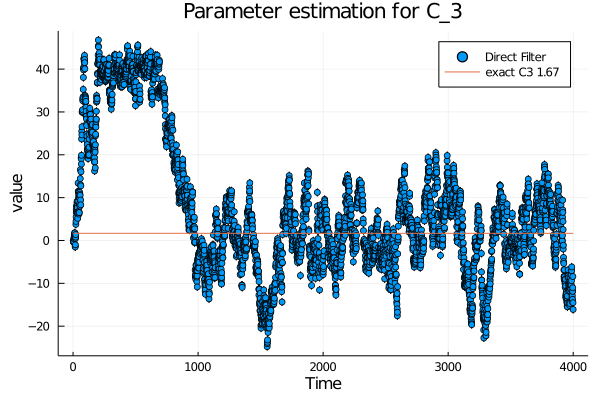}
  \caption{The estimation for paremeter $C_3$ }
  \label{fig:long}
  \end{figure}
 
 In the following experiment, we show that $C_3$ can be accurately estimated by the direct filter with more restrict environment setups. To this end, we let $\Sigma = diag(0.3^2, 0.8^2)$, and we observe the state process for $ 4000$ steps with the choice of depth as $D = 0.24$. The parameter estimation results are presented in Figure \eqref{fig:long}. We can see from this figure that by introducing a bigger portion of noise that allows exploration, the direct filter could eventually capture the true value of $C_3$ after approximately $1000$ data assimilation steps   
By dropping $1000$ burn-in step, the mean estimate of the direct filter for $C_3$ is $\tilde{C}_3 = 1.61$, which is fairly close to the true value of $C_3$,  i.e. $ C_3=1.671$. However, we also notice that it requires many more data assimilation steps to reach a good estimate, and the variance of the estimated parameter values is very large due to the large noise perturbation to the pseudo parameter process. This also indicates that $C_3$ is a less sensitive parameter in the tKdV model. Since our goal is to find the depth changes in shallow water, obtaining a good estimate for either $C_2$ or $C_3$ will give us a good understanding of $D$. Therefore, from the above experiments we know that it's better to use $C_2$ as an indicator for finding depth changes under our parameter estimation framework. 

%Part of the reason why it's harder to estimate $C_3$ than $C_2$ is that they vary at different scales: the scale of the noise required to be added to achieve equilibrium for the estimation of $C_3$ is much larger than $C_2$. 

\section{Numerical experiments for online detection of depth changes}
In this section, we apply the direct filter based parameter estimation to carry out online detection of depth changes. Since it's easier to estimate $C_2$ as we observed in the above section, in this practical application scenario we carry out parameter estimation for both $C_2$ and $C_3$ dynamically, but we only use the estimated results for $C_2$ to derive our estimate for $D$. As a result, our parameter estimation procedure will provide an \textit{online depth-change detection} method. The concept of ``online detection'' refers to the fact that the algorithm considers realtime online data one at a time for the filtering purpose. This would be a practical scenario for real life applications.

\subsection{One step change}
In the first example, we carry out a quick experiment to detect one abrupt depth change in the tank bottom topography, and we have the following setup. 
\begin{enumerate}
\item We assume that the exact signal process is generated by \eqref{rk4_data} perturbed by some noise.  
\item  Pick $\epsilon_n \sim \mathcal{N}(0,\Sigma)$, where $\Sigma=\text{diag}(0.3^2,0.017^2)$ 
	\item The upstream has $D^+=0.42$, so it is mildly shallow.
	\item The downstream has $D^+=0.24$, so it is shallower than the upstream. 
	\item We take the total time to be $1.2$ for both the upstream and downstream, and for data generation, we use RK4 with $\Delta t=0.0001$, so after the final time (T=2.4), we have collected $2.4\times10^4$ total datasets. %We use the terminal state before the ADC takes place as the initial state of the wave after ADC. 
	\item The window is picked to be containing $3500$ datasets. That is, after the burn-in period, we compute the moving average of every 3500 predictions as an estimate for $C_2$. 
\end{enumerate}

\vspace{0.5em}

\begin{figure}[hbt!]
  \centering
  \includegraphics[width=0.7\linewidth]{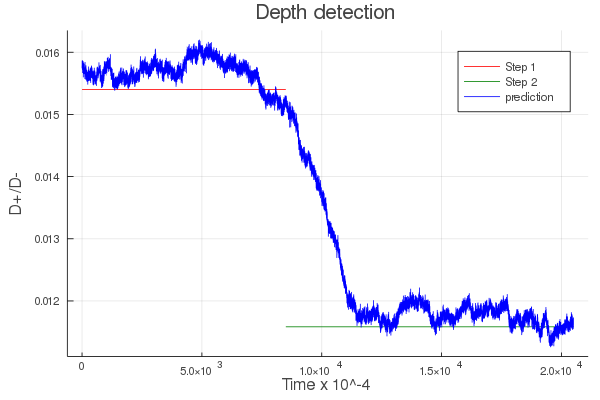}
  \caption{The estimation for paremeter $C_2$ }
  \label{fig:dd1}
  \end{figure}

The performance of parameter estimation by using the direct filter, which provides the depth estimation for $D$, is presented in Figure \ref{fig:dd1}.
In this experiment, the real abrupt depth change occurs at time $0.85$, where one can see the top red line changes to the bottom green line. The blue curve in the figure shows the online update of depth prediction -- based on the estimated values for $C_2$ discussed in Section \ref{DF-tKdV}. One can observe that the estimated water depth has been maintained at a good level until time $0.85$, when the curve plummets to the level of the green line. We can also see from this experiment that the direct filter immediately detected the depth change, and it could provide good estimates for the new depth fairly fast.

%It is also seen that since the window size we picked is relatively large, the burn in time after the ADC is also large. The prediction of $C^{down,pred}$ starts to reach equilibrium roughly at time $1.25$ ($1.6$ if we also consider the first 0.35 omitted), and it is observed to maintain at a good level. 

To study the accuracy of the method, we investigate the predicted values for $C_2$ in both the upstream and downstream region. For the upstream, we have true $C^{up}_2=0.154$, and the direct filter gives an estimated value $\tilde{C}^{up}_2=0.158$. 
The true downstream is $C^{down}_2=0.115$, and the direct filter estimated value is $\tilde{C}^{down}_2=0.117$. 
Hence, our direct filter method not only detected depth chang, but also gives good estimates of the depth level.

\subsection{Multiple depth change}
In this subsection, we carry out parameter estimation in the case that there are multiple depth changes. Recall the bottom topography presented in Figure \ref{fig:stepsplot}, it is observed that starting from the upstream, the water wave will go through various depth changes. We set up the problem as follows: 
\begin{enumerate}

	\item We assume that the exact signal process is generated by \eqref{rk4_data} perturbed by some noise. Starting from the left end, we have $D^0=1.0$ which means that we have a free surface and there is no step at the bottom. Then, we have $D^1=0.24, D^2=0.15, D^3=0.42$ according to the different bottom topography. Given an initial $\hat{u}_0$, we propagate it in the first region ($D^0=1.0$) and use the terminal state $\hat{u}_{T_1}$ as the initial state for the second region ($D^1=0.24$). Then, we take the terminal state in the second region as the initial state used for the third region ($D^2=0.15$). We repeat this same procedure until time runs out. 
	
	\item Take total time for each step region to be $1.2$, that is we let the wave travel in each region for a fixed time period of $T_i-T_{i-1}=1.2, i=1,2,3,4$. We take $\Delta t = 0.0001$, and we have carry out altogether $4.8 \times 10^4$ data assimilation steps. 
	\item Take $\epsilon_n \sim \mathcal{N}(0,\Sigma)$, where $(\Sigma)=\text{diag}(0.25^2,0.01^2)$ throughout the entire time horizon. 
\end{enumerate}  

\begin{figure}[hbt!]
  \centering
  \includegraphics[width=0.7\linewidth]{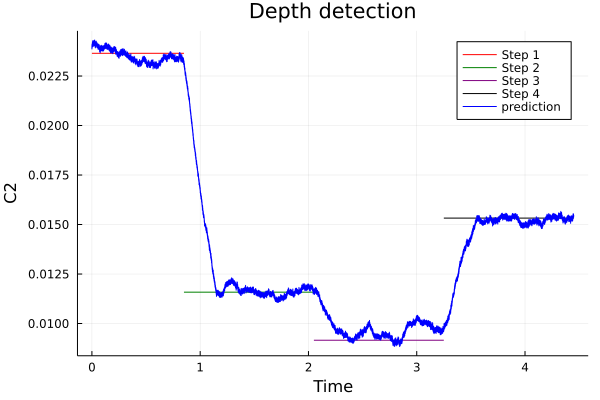}
  \caption{The estimation for paremeter $C_2$ }
  \label{fig:dd_multi}
  \end{figure}

By using the direct filter method, we present the depth estimation results in Figure \ref{fig:dd_multi},
where the horizontal line of different colors stand for the water depth. The size of the window picked for filtering is still $3500$. As one can observe, once an abrupt depth change occurs at the bottom topography, the parameter estimation results start to take a drastic change, and after the burn-in time, it will maintain at a relatively stable level after the moving  window completely shifts to the current region. The exact $C_2$ values for each steps (from left to right) are $0.236$, $0.0115$, $0.009$ and $0.0153$. And the the estimated values (after burn-in time) can give good approximations for the true values: $0.233$, $0.0113$, $0.00945$ and $0.0158$, respectively. 

Based on those estimated $C_2$ values, we compute by using \eqref{bigC} to obtain the estimated values for $D$:
\begin{equation}
	D= (C_2/ c_2)^2
\end{equation} 
In the following table, we list the true value for $D$, the estimated values and the error percentage.  
 \begin{table}\caption{Performance of the direct filter with different choice of $D$ }\label{Values:D:multi}
\centering%\begin{center}
 \begin{tabular}{||c c c||} 
 \hline
 $D$ & Estimated $D$ & Error $\%$ \\ [0.5ex] 
 \hline\hline
 1 & 0.978 & 2.2 \\ 
 \hline
0.24 & 0.228 & 5.0 \\
 \hline
 0.15 & 0.159 & 6.4 \\
 \hline
 0.42 & 0.446 & 6.1 \\ [1ex] 
 \hline
\end{tabular}
\end{table}
As a result, we see that the parameters related to the depth changes can be efficiently recovered by the direct filtering method.

\section{Conclusions}
In this work, we studied the ``inverse problem" discussed in \cite{BSM}, and we aim to estimate the water depth changes based on observations of the water surface displacements. To this end, we designed an optimal filtering based online parameter estimation technique by adopting the ``direct filter method" introduced in \cite{Bao_Direct}. We have demonstrated through numerical experiments that, by estimating the parameters in the tKdV model, we can effectively detect the water bottom topography. 

\bibliographystyle{apacite}
%\newpage

\section{Appendix}
This example is the first numerical example from \cite{Bao_Direct}. 

We consider the following discrete stochastic process: 
\begin{align}
	X_{n+1}^1 &= X_n^1 +(a_1 sin(X_t^2)+a_2 \frac{X^1_t}{1+|X^1_t|}) \Delta t +\sigma^1 \sqrt{\Delta t} W^1_t \\
	X_{n+1}^2 &= X_n^2 +(a_3 cos(X_t^1)+a_4 \frac{X^2_t}{1+|X^2_t|}) \Delta t +\sigma^2 \sqrt{\Delta t} W^2_t
\end{align} 
We pick $a_1=4.0$, $a_2=2.0$, $a_3=3.0$, $a_5=5.0$, $\sigma_1 = \sigma_2 =0.1$, $\Delta t =0.05$.

Define 
$$h: \mathbb{R}^2 \times \mathbb{R}^4 \rightarrow \mathbb{R}^2$$
\begin{equation}
h(X_n, \theta_{n+1})=
	\begin{pmatrix}
X_n^1 +(a_1 sin(X_t^2)+a_2 \frac{X^1_t}{1+|X^1_t|}) \Delta t\\
X_n^2 +(a_3 cos(X_t^1)+a_4 \frac{X^2_t}{1+|X^2_t|}) \Delta t
\end{pmatrix}
\end{equation}
where the input $\theta:=(a_1,a_2,a_3,a_4)$.

We assume the observed process takes the following form: 
$$ Y_{n+1}=H X_{n+1} + \eta_{n+1}=H h(X_n,\theta_{n+1}) + \sigma_Y \xi_{n+1}$$
Here $\sigma_Y$ is a diagonal matrix with each diagonal component $\sigma_Y^{(1)}$,$\sigma_Y^{(2)}$ determining the variance of the noise. $\eta$ is a random vector and from the above relationship, it is Gaussian with the following covariance matrix: 
\begin{equation}
	\eta^2 \sim diag((H_{11}\sigma^1  \sqrt{\Delta t})^2 + (\sigma_Y^{(1)})^2 ,(H_{22}\sigma^2 \sqrt{\Delta t})^2+ (\sigma_Y^{(2)})^2)
\end{equation}

We pick $H=\begin{pmatrix}
5& 0\\
0& 3
\end{pmatrix}$, $\sigma_Y=diag(0.1,0.1)$. 

The process for parameter process is assumed to take the following form: 
\begin{equation}
	\theta_{n+1}=\theta_n + \sigma_3*B_n
\end{equation}
where $\sigma_3$ is picked to be 0.1 and $B_n$ is a standard two dimensional Brownian Motion. 

By using the implementation step shown in the previous section, we obtain the estimated mean 
$$ \theta =\begin{pmatrix}
3.95\\
2.03\\
3.05\\
4.94
\end{pmatrix}$$
The following plot shows the estimated parameter process over 400 time steps. One can see that after a short amount of time, the process start to move close to some mean values which are very close to the exact values. 

\begin{figure}[h!]
\center
  \begin{subfigure}[b]{0.4\textwidth}
    \includegraphics[width=\textwidth]{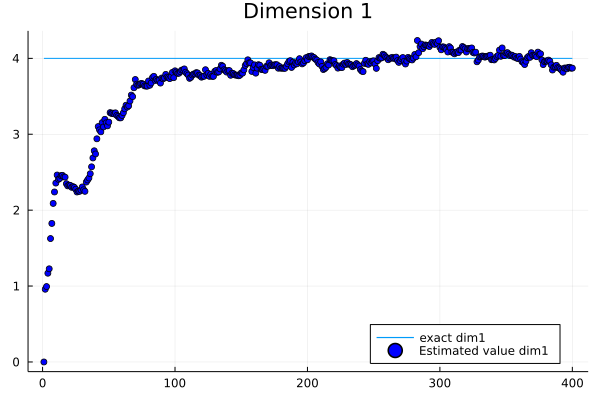}
    \caption{$a_1$}
    \label{fig:compare21}
  \end{subfigure}
  \begin{subfigure}[b]{0.4\textwidth}
  
    \includegraphics[width=\textwidth]{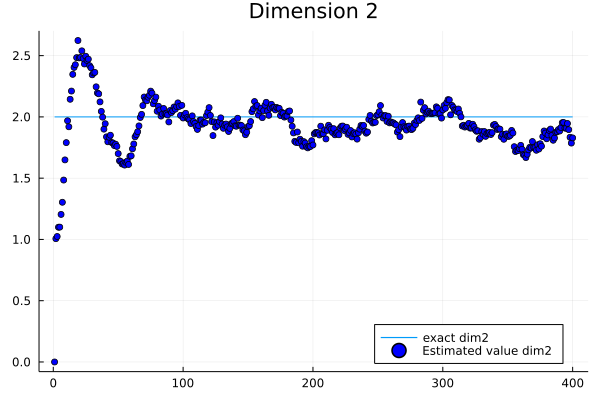}
    \caption{$a_2$}
    \label{fig:compare22}
  \end{subfigure}
\center
  \begin{subfigure}[b]{0.4\textwidth}
    \includegraphics[width=\textwidth]{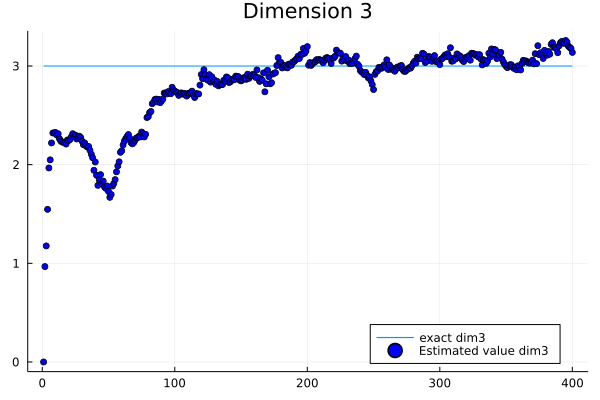}
    \caption{$a_3$}
    \label{fig:compare23}
  \end{subfigure}
  \begin{subfigure}[b]{0.4\textwidth}
    \includegraphics[width=\textwidth]{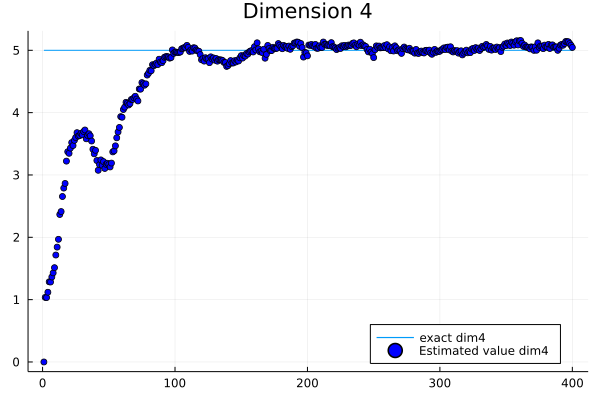}
    \caption{$a_4$}
    \label{fig:compare24}
  \end{subfigure}
  \caption{Parameter estimation for $a_1$,$a_2$,$a_3$ and $a_4$}
  \label{fig:compare}
\end{figure}

\end{document}